\documentclass[preprint,12pt]{elsarticle}
\usepackage{amsmath}
\usepackage{amssymb}
\usepackage{amsfonts}
\usepackage{lineno}

\usepackage{mathrsfs}

\usepackage[dvipdfm,colorlinks=true,linkcolor=blue,citecolor=blue,pdfstartview=FitH,pdftitle=title,pdfauthor=lixx]{hyperref}

\newtheorem{thm}{Theorem}

\newtheorem{definition}[thm]{Definition}
\newtheorem{prop}[thm]{Proposition}
\newtheorem{coro}[thm]{Corollary}

\newdefinition{rmk}{Remark}
\newdefinition{example}{Example}
\newproof{pf}{Proof}
\newproof{pot}{Proof of Theorem \ref{thm2}}

\journal{     }

\begin{document}
\begin{frontmatter}
\title{Moving Planes and Singular Points of Rational Parametric Surfaces}
\author[ustc]{Falai Chen\corref{cor1}}
\ead{chenfl@ustc.edu.cn}
\author[hfut,ustc]{Xuhui Wang}
\ead{wangxh05@mail.ustc.edu.cn}
\cortext[cor1]{Corresponding author} %
\address[ustc]{Department of Mathematics,
University of Science and Technology of China\\
Hefei, Anhui, 230026, China}
\address[hfut]{School of Mathematics,
Hefei University of Technology\\
Hefei, Anhui, 230009, China}

\begin{abstract}
In this paper we discuss the relationship between the moving planes
of a rational parametric surface and the singular points on it.
Firstly, the intersection multiplicity of several planar curves is
introduced. Then we derive an equivalent definition for the order of
a singular point on a rational parametric surface. Based on the new
definition of singularity orders, we derive the relationship between
the moving planes of a rational surface and  the order of singular
points. Especially, the relationship between the $\mu$-basis and the
order of a singular point is also discussed.
\end{abstract}

\end{frontmatter}


\section{Introduction}
Give an algebraic surface $f(x,y,z,w)=0$ in homogeneous form, the
singular points of the surface are the points at which all the
partial derivatives simultaneously vanish. Geometrically, A singular
point on the surface is a point where the tangent plane is not
uniquely defined, and it embodies geometric shape and topology
information of the surface. Detecting and computing singular points
has wide applications in Solid Modeling and Computer Aided Geometric
Design (CAGD).

To find the singular points of a parametric surface ${\mathbf
P}(s,t)$, one can solve ${\mathbf P}_s(s,t)\times {\mathbf
P}_t(s,t)={\mathbf 0}$. However, the find the orders of
singularities, one has to resort to implicit form. Let $f(x,y,z,w)$
be the implicit equation of ${\mathbf P}(s,t)$. A singular point
${\mathbf Q}=(x_0,y_0,z_0,w_0)$ of ${\mathbf P}(s,t)$ has order $r$
if all the partial derivatives of $f(x,y,z,w)$ with order up to
$r-1$ vanish at ${\mathbf Q}$, and at least one of the $r$-th
derivative of $f$ at ${\mathbf Q}$ is nonzero. Unfortunately,
converting the parametric form of a surface into implicit form is a
difficult task, which is still a hot topic of research
\cite{Buchberger, Buse1, Chen4, Chionh2, Cox-QSWB,Sed1}. Other
methods such as generalized resultants~ \cite{Buse2, Khetan2003,
Khetan2004} to find singular points don't need prior
implicitization. However, they are not designed for detecting and
computing the singular points of higher order (order $\geq 3$). In
this paper, we develop methods to treat singularities of rational
parametric surfaces directly. Specifically, we use moving surfaces
technique to study singularities of rational parametric surfaces and
the relationship between the order of the singularities and moving
surfaces.


The remainder of this paper is organized as follows. In Section 2,
we recalls some preliminary results about the $\mu$-basis of a
rational surface. The notion of intersection multiplicity of several
planar curves is also introduced. In Section 3, a new definition for
the order of a singular point directly from the parametric equation
of a surface is presented, and the equivalence of the new definition
with the classic definition is proved. In Section 4, we discuss the
relationship between the moving planes and $\mu$-basis of a rational
surface and the singular points of the rational surface. We conclude
this paper with future research problems.

\section{Preliminaries}
Let $\mathbb{R}[s,t]$ be the ring of bivariate polynomials in $s$,
$t$ over the set of real numbers $\mathbb{R}$. A rational parametric
surface in homogeneous form is defined
\begin{equation}\label{s-rsf}
{\mathbf P} (s,t)=\big( a(s,t),b(s,t),c(s,t),d(s,t) \big),
\end{equation}
where $a,b,c,d\in \mathbb{R}[s,t]$ are polynomials with
$\gcd(a,b,c,d)=1$. In order to apply the theory of algebraic
geometry, sometimes we need to work with homogeneous polynomials,
\begin{equation}\label{s-rsf-homo}
{\mathbf P} (s,t,u)=\big( a(s,t,u),b(s,t,u),c(s,t,u),d(s,t,u) \big).
\end{equation}

A \textit{base point} of a rational surface $\mathbf{P}(s,t)$ is a
parameter pair $(s_0, t_0)$ such that
$\textbf{P}(s_0,t_0)=\textbf{0}$.  Note that even if the rational
surface ${\mathbf P}(s,t)$ is real, the base points could be complex
numbers and possibly at infinity.

 A moving plane is a family of
planes with parametric pairs $(s,t)$~\cite{Sed1}
\begin{equation}\label{m-plane}
  A(s,t)x+B(s,t)y+C(s,t)z+D(s,t)
\end{equation}
where $A(s,t),B(s,t),C(s,t),D(s,t)\in \mathbb{R}[s,t]$. A moving
plane is said to \emph{follow} the rational surface \eqref{s-rsf} if
\begin{equation}
  A(s,t)a(s,t)+B(s,t)b(s,t)+C(s,t)c(s,t)+D(s,t)d(s,t)\equiv 0.
\end{equation}
The moving plane \eqref{m-plane} can be written as a vector form
$$\textbf{L}(s,t)=\big(A(s,t),B(s,t),C(s,t),D(s,t)\big). $$
Let $\textbf{L}_{st}$ be the set of the moving planes which follow
the rational surface $\textbf{P}(s,t)$, then $\textbf{L}_{st}$ is
exactly the syzygy module $syz(a,b,c,d)$ and is a free module of
rank $3$~ \cite{Chen4}.

A \emph{$\mu$-basis} of the rational surface \eqref{s-rsf} consists
of three moving planes $\mathbf{p}, \mathbf{q}, \mathbf{r}$
following \eqref{s-rsf} such that
$$[\mathbf{p}, \mathbf{q}, \mathbf{r}]=\kappa \mathbf{P}(s,t),$$
where $\kappa$ is nonzero constant and $[\mathbf{p}, \mathbf{q},
\mathbf{r}]$ is the outer product of $\mathbf{p}=(p_1,p_2,p_3,p_4)$,
$\mathbf{q}=(q_1,q_2,q_3,q_4)$, and $\mathbf{r}=(r_1,r_2,r_3,r_4)$
defined by
$$
[\mathbf{p},\mathbf{q},\mathbf{r}]=\Bigg(
\begin{vmatrix}
  p_2 & p_3 & p_4 \\
  q_2 & q_3 & q_4 \\
  r_2 & r_3 & r_4 \\
\end{vmatrix},
-\begin{vmatrix}
  p_1 & p_3 & p_4 \\
  q_1 & q_3 & q_4 \\
  r_1 & r_3 & r_4 \\
\end{vmatrix},
\begin{vmatrix}
  p_1 & p_2 & p_4 \\
  q_1 & q_2 & q_4 \\
  r_1 & r_2 & r_4 \\
\end{vmatrix},
-\begin{vmatrix}
  p_1 & p_2 & p_3 \\
  q_1 & q_2 & q_3 \\
  r_1 & r_2 & r_3 \\
\end{vmatrix}
\Bigg).
$$
A $\mu$-basis forms a basis for the syzygy module
$\textbf{L}_{st}$~\cite{Chen4}.

\begin{definition}\label{sing-def1}
Let $f(x,y,z,w)=0$ be the implicit equation of the parametric
surface ${\mathbf P}(s,t)$. Then $\mathbf{X}_0=(x_0,y_0,z_0,w_0)$ is
a singular point of {\em order} $r$ or an $r$-fold point if all the
derivatives of $f$ of order up to $r-1$ are zero at $\mathbf{X}_0$
and at least one $r$-th derivative of $f$ does not vanish at
$\mathbf{X}_0$. Specifically, ${\bf X}_0$ is a double point if and
only if
$$f_x({\bf X}_0)=f_y({\bf X}_0)=f_z({\bf X}_0)=f_w({\bf X}_0)=0,$$
and at least one of the second order derivatives is non-zero.
\end{definition}

To discuss the order of singular points on a rational parametric
surface, we need to recall some preliminary knowledge about the
intersection multiplicity of several curves in $\mathbb{P}^2
(\mathbb{C})$, the projective plane over the complex numbers.

\begin{definition}\label{curve-int-mul}\cite{Cox2001slide}
  Let $R$ be a local ring with  maximal ideal
  $\mathfrak{m}$ and $\mathcal{M}$ be a finitely generated
  $R$-module. Assume $R$ contains
  $k=R/\mathfrak{m}$. For $l\gg 0$, the Hilbert polynomial
  implies that
  $$
  \hbox{dim}_k (\mathcal{M}/m^{l+1} M)=\frac{e}{d!}l^d+\ldots,
  $$
  where $d=\hbox{dim}(R)$ and $e=e(\mathcal{M})$ is the
  multiplicity of $\mathcal{M}$.The refined case is as follows. Let $I$ be an ideal
  with $\mathfrak{m}^s \mathcal{M} \subset I\mathcal{M}$
  for some $s$, then $l\gg 0$ implies that
  $$
  \hbox{dim}_k (\mathcal{M}/I^{l+1} M)=\frac{\widetilde{e}}{d!}l^d+\ldots,
  $$
  where $\widetilde{e}=e(I,\mathcal{M})$ is the multiplicity of $I$ in $\mathcal{M}$.
\end{definition}

According to the above definition, we can define the intersection
multiplicity of several planar curves in $\mathbb{P}^2(\mathbb{C})$.
For planar curves $C_1, C_2, \ldots, C_v$ which are defined by
homogeneous equations $f_1(s,t,u)=0, \ldots, f_v(s,t,u)=0$
respectively. Let $\widetilde{f}_1, \ldots, \widetilde{f}_v $ be the
local equation of $C_1, C_2, \ldots, C_v$ near point $p$, then the
intersection multiplicity of these curves at point $p$ is
$$m(p)=e(I_p, R_p)$$
for $I_p=\langle \widetilde{f}_1, \widetilde{f}_2,\ldots,
\widetilde{f}_v \rangle$ and $R_p=\mathcal{O}_{\mathbb{P}^2,p}$,
which is the ring of rational functions defined at $p$.

Assume $\mathfrak{m}^s\subset I\subset R$, then
\begin{itemize}
  \item If $\mathfrak{m}^s\subset J \subset I \subset R$,
  then
  $e(J,R)\geq e(I,R)$.
  \item If $I^l J=I^{l+1}$, then
  $e(J,R)=e(I,R)$,  $J$ is a reduction ideal
  of $I$.
  \item If $I$ is generated by a regular sequence, then
  $e(I,R)=\dim_k R/I$,  $I$ is a complete
  intersection.
\end{itemize}
\begin{prop}\label{red-ideal}\cite{Cox2001slide}
\begin{enumerate}
  \item $I$ has a reduction ideal which is generated by a regular
  sequence.
  \item The regular sequence can be chosen to the generic linear
  combinations of the generators of $I$.
\end{enumerate}
\end{prop}

\section{The order of singular points on rational parametric surfaces}

Given a rational parametric surface, we first give a definition
about the order of singular points directly from the parametric
equation.

\medskip

\begin{definition}\label{sing-def2}
For a rational surface (\ref{s-rsf-homo}), a point
$\mathbf{X}_0=(x_0,y_0,z_0,w_0)$ (wlog, assume $w_0 \neq 0$) is a
$r$-fold singular point if
\begin{equation}\label{s-sg-def2}
\begin{aligned}
& w_0 a(s,t,u)- x_0 d(s,t,u)=w_0 b(s,t,u)- y_0 d(s,t,u)\\
& =w_0 c(s,t,u)- z_0 d(s,t,u)=0
\end{aligned}
\end{equation}
has $r+\lambda$ intersection points (counting multiplicity) in the
$(s,t,u)$ plane, where the multiplicity is defined in
\eqref{curve-int-mul}, and $\lambda$ is the number of base point of
the surface.
\end{definition}

We will show that the above definition is equivalent to the classic
definition of order of singularities (Definition~\ref{sing-def1}).

\begin{thm}\label{c2-sg-same}
Definition \ref{sing-def1} and Definition \ref{sing-def2} are
equivalent.
\end{thm}
\begin{pf}
$\mathbf{X}_0$ is an $r$-fold singular point on a surface if and
only if, for a generic line passing through $\mathbf{X}_0$, the line
intersects the surface at $n-r$ distinct points besides
$\mathbf{X}_0$, here $n$ is the implicit degree of the surface.
Without loss of generality, assume the singular point is at the
origin $\mathbf{X}_0=
(x_0,y_0,z_0,w_0)=(0,0,0,1)$. Let the generic line be defined by %
$$l=L_1 \cap L_2,$$
where
\begin{align*}
  &L_1: \ \alpha_1 x +\alpha_2 y+ \alpha_3 z=0,&\\
  &L_2: \ \beta_1 x +\beta_2 y+ \beta_3 z=0.&
\end{align*}
Consider the two planar curves
\begin{align*}
  &C: g_1=\alpha_1 a(s,t,u)+\ldots +\alpha_3 c(s,t,u)=0,&\\
  &D: g_2=\beta_1 a(s,t,u)+\ldots +\beta_3 c(s,t,u)=0.&
\end{align*}
Let $Z$ be the common zeros of $a,b,c$, $S$ be the surface, $l^{-}$
be the line segments by removing the origin from the line $l$,
denote $\varphi: (s,t,u)\mapsto \mathbf{P}(s,t,u)$, then
$$C\cap D= \varphi^{-1}(S \cap l^{-}) \cup Z.$$

By Bezout's theorem, one has
\begin{equation}\label{c2-bezout}
n^2=\#(S \cap l^{-})+\sum_{p\in Z} \dim \mathcal{O}_{p}/\langle g_1,
g_2\rangle_{p}.
\end{equation}

From \eqref{red-ideal}, we can get that ${g_1}_p$ and ${g_2}_p$ is a
reduction ideal of $\langle a(s,t),b(s,t), c(s,t)\rangle_p$, where
$p$ is the intersection point of $g_1=0,g_2=0$. Thus,
$$e(\langle a,b,c\rangle_{p}, \mathcal{O}_{p})=
e(\langle g_1, g_2\rangle_{p}, \mathcal{O}_{p})=\dim_{k} \mathcal{O}_{p}/\langle {g_1}, {g_2} \rangle_p$$
Therefore, (\ref{c2-bezout}) is equivalent to
\begin{equation}\label{c2-bezout2}
n^2=\#(S \cap l^{-})+\sum_{p\in Z} e(\langle a,b,c\rangle_{p},
\mathcal{O}_{p}).
\end{equation}

Let $Z_1$ be the point set which satisfies $a=b=c=0$ and $d\neq0$,
$Z_2$ be the point set which satisfies $a=b=c=d=0$, then
$$Z=Z_1 \cup Z_2, \ \ \hbox{and }\ Z_1 \cap Z_2=\emptyset.$$
Therefore, (\ref{c2-bezout2}) is equivalent to
\begin{equation}\label{c2-bezout3}
n^2=\#(S \cap l^{-})+r+\lambda.
\end{equation}
Since the implicit degree of the surface is $n^2-\lambda$, we
immediately get that Definition~\ref{sing-def1} and Definition
~\ref{sing-def2} are equivalent. \qed
\end{pf}

\section{Relationship between moving planes  and  singular points}

In this section, we study the relationship between the moving planes
and the order of singular points on a rational parametric surface.

\begin{thm}\label{c2-sg-mp}
Let $\mathbf{P}(s,t,u)$ be a parametric surface with no base points,
and $\mathbf{L}(s,t,u)$ be a moving plane following
$\mathbf{P}(s,t,u)$. If $\mathbf{X}_0=(x_0,y_0,z_0,w_0)$ (assume
$w_0\neq 0$) is an $r$-fold singular point on the surface, then
\begin{align*}
&w_0 a(s,t,u)-x_0 d(s,t,u)=0,\quad w_0 b(s,t,u)-y_0
d(s,t,u)=0,&\\&w_0 c(s,t,u)-z_0 d(s,t,u)=0,\quad
\mathbf{L}(s,t,u)\cdot \mathbf{X}_0=0&
\end{align*}
 have $r$ intersection points (counting multiplicity).
\end{thm}
\begin{pf}
The rational parametric surface $\mathbf{P}(s,t)$ has the following
special three moving planes
\begin{align*}
  &\mathbf{L}_1:=\big(-d(s,t),0,0,a(s,t)\big),&\\
  &\mathbf{L}_2:=\big(0,-d(s,t),0,b(s,t)\big),&\\
  &\mathbf{L}_3:=\big(0,0,-d(s,t),c(s,t)\big),&
\end{align*}
and they belong to $\mathbf{L}_{s,t}$. Given a moving plane
$\mathbf{L}(s,t)=\big(A(s,t),B(s,t),C(s,t)$ $D(s,t)\big)$ follow the
rational surface $\mathbf{P}(s,t)$, assume
$A(s,t),B(s,t),C(s,t),D(s,t)$ are relatively prime (If they are not
relatively prime, we can deal with it similarly). As four
dimensional vectors, $\mathbf{L}_1, \mathbf{L}_2, \mathbf{L}_3$ are
all perpendicular to $\mathbf{P}(s,t)$, and $\mathbf{L}$ is also
perpendicular to $\mathbf{P}(s,t)$. Thus, there exist $h, h_1, h_2,
h_3\in \mathbb{R}[s,t]$, and $\gcd(h_1,h_2, h_3)=1$, such that
\begin{align}\label{c2-relation}
h \mathbf{L}(s,t)&=h_1 \mathbf{L}_1(s,t) +h_2 \mathbf{L}_2(s,t)
+h_3  \mathbf{L}_3 (s,t)&\\
&= (-h_1 d,-h_2 d, -h_3 d, h_1 a+h_2 b +h_3 c).\nonumber&
\end{align}
Since $\gcd(A,B,C,D)=1$ and $\gcd(h_1, h_2, h_3)=1$,
$$
h=\gcd(-h_1 d,-h_2 d, -h_3 d, h_1 a +h_2 b +h_3 c)=\gcd(d, h_1 a+h_2
b +h_3 c).
$$
Thus $h|d$.

For an $r$-fold singular point $\mathbf{X}_0$ on the surface,
\eqref{s-sg-def2} have $r$ intersection points (counting
multiplicity):
\begin{equation}\label{c2-int-mul}
(p_0, \ldots, p_v)=\big((s_0,t_0,u_0),\cdots, (s_v,t_v,u_v)\big),
\end{equation}
and the multiplicity at $p_i$ is $m_i , i=0,\cdots, v$, thus
$r=m_0+\ldots+m_v$. From \eqref{c2-relation}, we have
\begin{equation}\label{mp-smp}
h(s,t) l(s,t)=h_1(s,t) l_1 (s,t)+h_2(s,t) l_2(s,t)+h_3(s,t)
l_3(s,t),
\end{equation}
where $l(s,t)=\mathbf{L}\cdot \mathbf{X}_0, l_i(s,t)=\mathbf{L}_i
\cdot \mathbf{X}_0,i=1,2,3$.

\bigskip

Now we discuss the two possible cases.\\
$\bullet$ If the intersection point $p_i$ is not at infinity, assume
$p_i=(s_i,t_i,1)$, and from \eqref{mp-smp}, we have
\begin{align*}
&h(s-s_i,t-t_i) l(s-s_i,t-t_i)=h_1(s-s_i,t-t_i) l_1
(s-s_i,t-t_i)&\\
&+h_2(s-s_i,t-t_i) l_2(s-s_i,t-t_i)+h_3(s-s_i,t-t_i)
l_3(s-s_i,t-t_i).&
\end{align*}
Thus,the ideals generated by local equations of $h(s,t,u) l(s,t,u),$
$l_1(s,t,u),$ $l_2(s,t,u), l_3(s,t,u)$ near $p_i$ and local
equations of $l_1 (s,t,u),$ $l_2(s,t,u),$ $l_3(s,t,u)$ near $p_i$
are same.\\
$\bullet$  If  intersection point $p_i$ is at infinity, assume
$p_i=(1,t_i,0)$. Since $\mathbf{L}(s,t)\cdot \mathbf{P}(s,t)\equiv
0$, its homogeneous form also satisfies:
$$
\mathbf{L}(s,t,u) \cdot \mathbf{P}(s,t,u)\equiv 0,
$$
and the dehomogenized form also has $ \mathbf{L}(1,t,u) \cdot
\mathbf{P}(1,t,u)\equiv 0.$ Similar to the above analysis,  their
exist $h',h_1',h_2',h_3'\in \mathbb{R}[t,u]$, and
$\gcd(h_1',h_2',h_3')=1$, such that
$$
h'~l(1,t,u)=h_1' l_1 (1,t,u)+h_2' l_2(1,t,u)+h_3' l_3(1,t,u),
$$
where $l(s,t,u),l_i(s,t,u)$ is the homogeneous form of
$l(s,t),l_i(s,t), i=1,2,3$, and $h'|d(1,t,u)$. Therefore, the ideals
generated by local equations of $h' l(1,t,u),$ $l_1 (1,t,u),$
$l_2(1,t,u),l_3(1,t,u)$ near $p_i$ and local equations of $l_1
(1,t,u),$ $l_2(1,t,u),$ $l_3(1,t,u)$ near $p_i$ are also same.

\bigskip

Since $h(p_i)\neq 0, h'(p_i)\neq 0$ (otherwise, $p_i$ is a base
point of $\mathbf{P}(s,t)$), and based on the definition of the
multiplicity,  $p_i,i=1,\ldots,v$ are also the intersection points
of
\begin{align*}
&w_0 a(s,t,u)-x_0 d(s,t,u)=0,\quad w_0 b(s,t,u)-y_0
d(s,t,u)=0,&\\&w_0 c(s,t,u)-z_0 d(s,t,u)=0,\quad
\mathbf{L}(s,t,u)\cdot \mathbf{X}_0=0.&
\end{align*}
and the intersection multiplicity at $p_i$ are also same. \qed
\end{pf}

{\bf Remark} 1. $\mu$-basis $\mathbf{p}, \mathbf{q}, \mathbf{r}$ of
$\mathbf{P}(s,t)$ are also three special moving planes of the
surface, and from Theorem \ref{c2-sg-mp},
\begin{equation}
\begin{split}
  &\mathbf{L}_1 (s,t,u)\cdot \mathbf{X}_0 =0, \  \mathbf{L}_2 (s,t,u)\cdot
\mathbf{X}_0 =0,\   \mathbf{L}_3 (s,t,u)\cdot \mathbf{X}_0 =0, \\
& \mathbf{p} (s,t,u)\cdot \mathbf{X}_0 =0,\  \   \mathbf{q}
(s,t,u)\cdot \mathbf{X}_0 =0,\  \ \ \mathbf{r} (s,t,u)\cdot
\mathbf{X}_0=0
\end{split}
\end{equation}
also have the $r$ intersection points.

\bigskip

Next we discuss the relationship of the $\mu$-basis and the order of
singular points on a rational parametric surface.

For any moving plane $\mathbf{l}(s, t) \in \mathbf{L}_{s,t}$, there
exist polynomials $h_i (s, t), i = 1, 2, 3$, such that
$\mathbf{l}(s, t) = h_1 \mathbf{p} + h_ 2 \mathbf{q} + h_3
\mathbf{r}$ (Theorem 3.2 in \cite{Chen4}), thus
$$
\langle \mathbf{L}_1\cdot \mathbf{X}_0, \mathbf{L}_2\cdot
\mathbf{X}_0, \mathbf{L}_3\cdot \mathbf{X}_0, \mathbf{p}\cdot
\mathbf{X}_0, \mathbf{q}\cdot \mathbf{X}_0, \mathbf{r}\cdot
\mathbf{X}_0\rangle = \langle \mathbf{p}\cdot \mathbf{X}_0,
\mathbf{q}\cdot \mathbf{X}_0, \mathbf{r}\cdot \mathbf{X}_0 \rangle.
$$
Similar to the proof of Theorem \ref{c2-sg-mp}, we can get that
\begin{equation}
  \mathbf{p}(s,t,u)\cdot \mathbf{X}_0=\mathbf{q}(s,t,u)\cdot
\mathbf{X}_0=\mathbf{r}(s,t,u)\cdot \mathbf{X}_0=0
\end{equation}
also have the same multiplicity at each intersection point $p_i,
i=0,\ldots, v$ as Equation \eqref{s-sg-def2}.

\bigskip

An immediate consequence of the above analysis is:
\begin{coro}\label{c2-sg-order-mu}
For a rational parametric surface $\mathbf{P}(s,t)$ with no base
point and its $\mu$-basis $\mathbf{p}, \mathbf{q}, \mathbf{r}$, then
$\mathbf{X}_0=(x_0,y_0,z_0,w_0)$ is $r$-fold singular point on the
surface if and only if the intersection points of
$$\mathbf{p}(s,t,u)\cdot \mathbf{X}_0= \mathbf{q}(s,t,u)\cdot \mathbf{X}_0= \mathbf{r}(s,t,u) \cdot \mathbf{X}_0=0$$
is $r$ (counting multiplicity).
\end{coro}

\bigskip

Now we consider the moving surface ideal:
$$
I'=\langle dx-a,dy-b,dz-c,dw-1 \rangle \cap \mathbb{R}[x,y,z,s,t].
$$
Theorem 3.4, 3.5 of \cite{Chen4} shows that $I'$ is a prime ideal
and $f(x,y,z,s,t)\in I'$ if and only if $f(x,y,z,s,t)=0$ is a moving
surface following the rational surface $\mathbf{P}(s,t)$. Moreover,
if $\mathbf{P}(s,t)$ contains no base point then
\begin{equation}
  I'=\langle p, q, r\rangle
\end{equation}
where $p=\mathbf{p}\cdot (x,y,z,1)$, $q=\mathbf{q}\cdot (x,y,z,1)$,
$r=\mathbf{r}\cdot (x,y,z,1)$. Therefore, For a rational surface
$\mathbf{P}(s,t)$ with no base point, and any moving surface $f=0$
following it, $f\in \langle p, q, r\rangle$. Based on the above
analysis, we can improve Theorem 6 as the following theorem
\begin{thm}\label{c2-sg-mq}
For a parametric surface $\mathbf{P}(s,t,u)$ with no base point and
any moving surface $f(x,y,z,w,s,t,u)=0$ following it. If
$\mathbf{X}_0=(x_0,y_0,z_0,w_0)$ (assume $w_0\neq 0$) is an $r$-fold
singular point on the surface, then
\begin{align*}
&w_0 a(s,t,u)-x_0 d(s,t,u)=w_0 b(s,t,u)-y_0 d(s,t,u)&\\&=w_0
c(s,t,u)-z_0 d(s,t,u)=f(x_0,y_0,z_0,w_0,s,t,u)=0&
\end{align*}
 have and only have $r$ intersection points (counting multiplicity), where $f(x,y,z,$ $w,s,t,u)$ is the homogeneous form of $f(x,y,z,$ $w,s,t)$.
\end{thm}

\section{Conclusion}
To make the $\mu$-basis more applicable in computing the singular
points, we discuss the relations between moving planes (Specially,
$\mu$-basis) and singular point of the rational surface.

In the future, we will discuss how to detect and compute singular
points on an rational parametric surface based on moving planes (or
$\mu$-basis) in an efficient way.

\section{Acknowledgements}
This work was partially supported by the NSF of China grant
10671192.

\end{document}